\documentclass[opre,nonblindrev,copyedit]{informs2}

\usepackage{endnotes}
\let\footnote=\endnote

\TheoremsNumberedThrough     
\ECRepeatTheorems

\EquationsNumberedThrough    

\usepackage{mathtools}

\usepackage{appendix}

\newcommand{\pr}{\mbox{\sf P}}
\newcommand{\ex}{{\bf\sf E}}               
\newcommand{\var}{\mbox{\sf Var}}

\newcommand{\al}{\alpha}                
\newcommand{\sig}{\sigma}               
\newcommand{\s}{\sigma}               

\newcommand{\startb}{\parindent0pt\bf}  
\newtheorem{thm}{Theorem}

\newtheorem{pro}[thm]{Proposition}

\newtheorem{defn}[thm]{Definition}

\def\eqd{{\buildrel {\rm d} \over =}}

\def\eqd{{\buildrel\rm d\over =}}

\def\hf{\frac{1}{2}}

\begin{document}

\RUNAUTHOR{Lu, Squillante \& Yao}

\RUNTITLE{Matching Supply and Demand}

\TITLE{Matching Supply and Demand in Production-Inventory
Systems: Asymptotics and Optimization}

\ARTICLEAUTHORS{%
\AUTHOR{Yingdong Lu, Mark S.\ Squillante}
\AFF{Mathematical Sciences Department, IBM Research, Yorktown
  Heights, NY 10598, USA.  \EMAIL{\{yingdong,mss\}@us.ibm.com}} 
    \AUTHOR{David D. Yao}
    \AFF{Department of Industrial Engineering and Operations Research, Columbia University,
    New York, NY 10027, USA.  \EMAIL{yao@ieor.columbia.edu}}
} 

\ABSTRACT{%
We consider a general class of high-volume, fast-moving production-inventory
systems based on both lost-sales and backorder inventory models.
Such systems require a fundamental understanding of the asymptotic behavior
of key performance measures under various supply strategies, as well as the
pre-planning of these strategies. Our analysis relies on a thorough study of
the asymptotic behavior of a random walk with power drift,
which is of independent interest.
In addition to providing key insights, our analysis leads to approximations of
the corresponding optimization problem that yield simple solutions which are
close to optimal.
We also establish an equivalence between the lost-sales and backorder models
when both have the same penalty cost that becomes large.
}




\maketitle


\baselineskip 18pt

\section{Introduction}
\label{sec:intro}

The production-inventory problem we study is that of matching supply and demand
over a planning horizon of $N$ periods, with the provision that a set of supply
decisions must be made before the first period.
For example, the production-inventory planning horizon might represent a
``selling season'' and the decisions to match supply with demand are made before
the season starts.
This is not merely a quantity decision --- how much
to stock up so as to fulfill demand over the  entire season; rather,
it is more like a capacity planning decision --- how to set the
supply or production rates as a function of time  over the horizon, vis-\`a-vis
the demand forecast.
We consider this production-inventory matching problem within the context of both
lost-sales and backorder inventory models.

More precisely, consider the following general supply and demand matching problem.
Suppose $D_1+\cdots +D_n$ denotes the cumulative demand up to period
$n=1,2,\dots, N$; and assume $\var (D_n)=\sig^2$ for all $n$.
Further suppose the set of supply decisions before the first period takes
the form $\ex[D_1]+\cdots + \ex[D_n]+\kappa\sig n^\al$, where  $\al$ and
$\kappa$ are decision variables such that $\al\in [0,1]$ while $\kappa$
can be either positive or negative.
This assumed form allows a fair amount of flexibility:
one can choose to over-supply (positive $\kappa$), thus
creating safety stock, if the lost-sales penalty is high;
or to under-supply (negative $\kappa$), if the production cost
or holding cost is high.
We are especially interested in asymptotic results.
When $n$ is large, demand will approach a normal distribution, with
a standard deviation equal to $\sig\sqrt{n}$.
Shall we match supply with demand by setting $\al=\hf$?
What is the corresponding lost-sales quantity?
What is the impact on lost sales if $\al$ takes on other values?
What is the right value for $\kappa$;
in particular, should it be positive or negative?

We note that this general supply and demand matching problem arises in a wide
variety of production-inventory systems.
As such, there is another way to interpret the asymptotics associated with
large values of $N$ in our model.
Instead of viewing this as a prolonged planning horizon, we can view the latter
as having a fixed length, but with the demand being scaled up.
In other words, the original total demand for the entire season, $D$, now
becomes $D_1+\cdots +D_N$, with the $D_n$ being independent copies of $D$.
The capacity problem continues to consist of making the set of supply decisions
for the entire horizon in advance, with the corresponding deliveries received
every period to satisfy the demand.
Indeed, this alternative view is closer to the class of production-inventory
applications we have in mind --- a high-volume, fast-moving supply and demand
context, which requires an understanding of the asymptotic behavior of key
system performance measures under various supply strategies, as well the
pre-planning of such strategies.

As a specific business analytics instance of such capacity problems often encountered within
the context of workforce management for large services providers, we note
that any significant workforce actions aimed at altering a planned resource
supply usually require long lead times, whereas customer demand can exhibit
volatility and uncertainty over relatively short periods of time.
Competition within the marketplace and the dynamics within these companies
create constant pressure for growth; under this pressure, workforce capacity
planning routinely examines resource requirements when the demand forecast
is scaled up in different proportions.
These considerations in part motivate the current problem setting.
Another application of a similar flavor concerns capacity provisioning in
cloud computing environments where high-frequency demand having considerable
uncertainty needs to be fulfilled with a predetermined supply of resources.
Any changes in these supply decisions are costly and should be avoided.

Now, let $S_n$ denote the net demand at period $n=1,\ldots,N$, namely the
difference between demand and supply (both in terms of the cumulative quantities)
up to period $n$.
In the backorder model, where unmet demand in any period will be backlogged
and supplied by future supplies, backorder and inventory are, respectively,
the positive and negative parts of $S_n$; as such, the model is quite tractable.
On the other hand, the lost-sales model, in which unmet demand in any period will be lost,
is much harder to solve because the cumulative lost-sales quantity, $L_n$, is the running
maximum of $S_n$: $L_n=\max_{0\le j\le n}\{S_j\}$ (with $S_0:=0$), which depends on $S_1, S_2,\ldots, S_n$.
When supply is linear, i.e., $\al=1$, then $S_n$ is a random walk and $L_n$ is
a well-studied object; in particular, $\ex[L_n]$ is well-known and very
accessible.
For other values of $\al$, however, $\ex[L_n]$ does not appear to have explicit
expressions.
The only possible exception is $\al=\hf$, in which case one can approximate
$S_n$ by a Brownian motion with square-root drift.
We derive such an analysis by relating this Brownian motion to the hitting
time of an Ornstein-Uhlenbeck process, for which the density function can be
obtained by inverting Laplace transforms.

Given this level of difficulty, our approach consists of constructing upper-
and lower-bounds for $\ex[L_N]$, and analyzing the asymptotic behavior of
these bounds for large $N$.
Our findings for such asymptotic behaviors are summarized as follows.
When $\kappa >0$ (the over-supply case), setting $\al$ within $(\hf,1]$ will
result in a bounded $\ex[L_N]$, whereas $\al \in [0,\hf]$ will cause $\ex[L_N]$
to grow on the order of $\sqrt {N}$.
When $\kappa < 0$ (the under-supply case), $\ex[L_N]$ will grow on the order of
$N^\al$ upon setting $\al$ within $(\hf,1]$; whereas it will grow on the order
of $\sqrt {N}$ when $\al \in [0,\hf]$.

In terms of minimizing the total cost over both production cost and lost-sales
penalty, our results imply that setting $\al=\hf$ is the only meaningful choice
leading to non-trivial results, independent of the values of the penalty and holding costs.
Hence, what remains is to find the best value of $\kappa$.
We first obtain an asymptotically optimal solution through a
Brownian approximation of the objective function for large $N$.
Then, for each finite $N$, we examine optimal solutions through the upper
and lower bounds of $\ex[L_N]$ as surrogates for the original problem.
Both bounds lead to simple optimal solutions, and our results include performance
guarantees for both problems as approximations to the original problem.
We also show that these results are readily extended to include holding costs.
In addition, we establish the equivalence between the lost-sales model and
the backorder model when both have the same penalty cost that goes to infinity.

\subsection{Related Work}
A brief review of the related literature is in order.
To the extent that our set of supply decisions has to be made before the
season starts, our model resembles the newsvendor model.
There are important distinctions, however.
To start with, the newsvendor decision is a stocking quantity decision, whereas
we are concerned with the supply {\it dynamics} over time --- all $N$ periods
of the planning horizon.
Consequently, the lost-sales quantity in the newsvendor model is a single number that
can be accounted for at the end of the horizon.
In our model, lost sales build up over the horizon, as demand and supply are
realized and evolve.
Granted, under our cost model which minimizes the total production cost and
lost-sales penalty over the horizon, the lost-sales quantity $\ex[L_N]$ is
associated with the end of the horizon, as in the newsvendor model.
However, we consider inventory holding costs as well the backorder alternative
to lost sales, features that are not present in the newsvendor model.

The lost sales problem is known to have a much more complicated structure than
the problem of backorder, and hence it is more difficult to analyze.
Regarding the optimal policy, the exact description is not known even in the
simplest case.
Karlin and Scarf~\cite{KarlinScarf} demonstrate that the base-stock policies
are not optimal even for systems with a lead time of one period.
Monotonicity properties of the optimal ordering policy are derived by
Morton~\cite{Morton}, Zipkin~\cite{zipkin} and most recently
Huh~\cite{HuhOR10}, each through different methods.
Huh {\it et al.}~\cite{HuhMS09} show that, asymptotically as the ratio of
unit penalty cost and unit holding cost grows to infinity, the difference
in performance between a particular base-stock policy and the optimal policy
will eventually vanish.
Thus, base-stock policies can be regarded as asymptotically optimal when there
is such a significant imbalance between the penalty cost and holding cost.

Levi {\it et al.}~\cite{levi} study production and lost sales over multiple
periods focusing on the performance of a simple heuristic policy that balances
the inventory holding cost and lost-sales penalty in each period, with the aim
of establishing its robustness relative to the optimal policy in a very general
context.
The setting of Levi {\it et al.} is dynamic programming: in each period a
different decision can be made and actions taken accordingly.
In contrast, we allow only a single set of decisions at the beginning of the
horizon;
and all subsequent costs due to production, lost-sales or backlog, and holding
inventory are consequences of this set of decisions.
The way we quantify our heuristic solution (from the upper-bound problem),
however, appears to have a similar flavor to the proven ``2-approximation''
status of the dual-balancing rule in \cite{levi}.
Specifically, the performance of our heuristic, via the ratio of the upper-
and lower-bound solutions, is ``2-plus'' --- where this slightly worse
performance guarantee than a 2-approximation algorithm can be attributed to the
fact that there is no recourse over the planning horizon in our model, which,
to make matters worse, can be infinitely long.
A complementary set of results for the lost-sales model has been recently
obtained by Goldberg {\it et al.}~\cite{goldberg}.

Another body of research related to our analysis concerns process flexibility;
see the recent paper of Chou {\it et al.}~\cite{chou} and the references
therein.
Specifically, process flexibility can be modeled by a random walk with
functional drift, called a generalized random walk in~\cite{chou}.
When proper functional forms are taken, it can be seen that the cost structure
under full flexibility has the same form as the cost structure of our backorder
model, while the cost structure under a so-called chaining option corresponds
to the structure of our lost-sales model.
Upon applying the methodology developed in this paper, the effectiveness of
chaining can be demonstrated and quantified for more general processes.

\subsection{Paper Organization}
The remainder of the paper is organized as follows.
In Section~\ref{sec:prelim}, we start with a formulation of our models ---
while focusing on the lost-sales model, we also cover the backorder model,
followed by a derivation of the main performance measures associated with
these models.
We then present preliminaries regarding the normal distribution function
and the loss function, and derive estimates for $\ex[L_N]$ in the case of
$\al=1$.
In Section~\ref{sec:bounds}, we continue with estimates of $\ex[L_N]$ for all
other cases, through constructing upper- and lower-bounds and also
deriving a corresponding Brownian approximation for the case of $\al=1/2$.
Optimization models that minimize total cost over production, lost sales
(or backorder) and inventory are studied in Section~\ref{sec:optimality},
followed by concluding remarks in Section~\ref{sec:conclusions}.


\section{Problem Formulation and Preliminaries}
\label{sec:prelim}

\subsection{Key Performance Measures}

Consider a planning horizon that consists of $N$ periods,
indexed by $n=1,\dots, N$.
For each period $n$, let $D_n$ denote the demand, a random variable;
and $x_n$ the production or supply quantity, a decision variable.
At the beginning of the planning horizon, we need to determine the
production quantities for all $N$ periods, as part of the capacity
planning problem.
Once this set of decisions is made, then in each period $n$ there will
be $x_n$ units available along with any units left over from previous
periods to supply demand $D_n$.
In the lost-sales model,
any demand surpassing the available supply is lost by the end of each period;
whereas in the backorder model, any unfilled demand will be
backlogged, to be supplied in a later period.
Correspondingly, the quantities of interest include the following.
\begin{itemize}
\item $L_n$: the cumulative demand
shortfall, in the lost-sales model, up to period $n$;
\item
$B_n$: the backlogged demand, in the backorder model,  at the end of each
period $n$;
\item
$H_n$ and $H'_n$: the on-hand inventory left at the end of each period $n$,
respectively for the lost-sales and backorder models.
\end{itemize}

\begin{pro}
\label{pro:lost}
{\rm
Define $S_0:=0$ and $S_n:=D_1-x_1+\cdots +D_n-x_n$.
Further define $x^+:=\max\{x,0\}$ and $x^-:=-\min\{x,0\}$.
Then, for $n=1,\dots, N$, we have
\begin{eqnarray}
B_n= S_n^-, &\qquad& H'_n=S_n^+, \label{back}\\
L_n=\max_{0\le j\le n}\{S_j\}, &\qquad&  H_n = L_n-S_n.
\label{lost}
\end{eqnarray}
}
\end{pro}

{\startb Proof.}
In the backorder model, any unmet  demand (if $S_n>0$) or any leftover
inventory (if $S_n<0$) at the end of each period $n$ will be carried over;
hence, leading to the expressions in (\ref{back}).

The $H_n$ expression in (\ref{lost}) follows from
\begin{eqnarray*}
H_n=\sum_{i=1}^n x_i -\left(\sum_{i=1}^n D_i-L_n\right),
\end{eqnarray*}
i.e., total production minus total {\it supplied} demand equals inventory.

For the $L_n$ expression  in (\ref{lost}), note that the result
holds trivially for $n=1$: $L_1=\max\{0,D_1-x_1\}$.
Suppose it holds for $n$. Then,
we have $L_n=S_k$ for some $k\le n$,
which means that there will be no more shortage for periods $k+1$ through $n$.
Therefore, the units available to supply $D_{n+1}$ consist of $x_{n+1}$ and
the leftover from periods $k+1$ through $n$, namely
$$-(D_{k+1}-x_{k+1} + \cdots +D_n-x_n).$$
(Note that the above expression is non-negative, as the quantity in parentheses
cannot exceed $0$ due to no shortage from periods $k+1$ through $n$.)
Hence, at the end of period $n+1$, shortage will occur if
$$D_{n+1}- [x_{n+1} - (D_{k+1}-x_{k+1} + \cdots +D_n-x_n)] >0 ,$$
in which case the cumulative shortage will be adding $S_k$ to the above expression,
with the sum being equal to $S_{n+1}$.
Otherwise, i.e., if the left-hand side above is $\le 0$, then
the cumulative shortage stays at $S_k$.
In summary, we have
\begin{eqnarray*}
L_{n+1}= \max\{S_k,S_{n+1}\}= \max\{L_n, S_{n+1}\}= \max_{0\le j\le n+1}\{S_j\},
\end{eqnarray*}
which is as desired.
 \hfill$\Box$

\medskip

In what follows, we shall focus on $L_n$, since $H_n$ relates directly to $L_n$;
whereas $B_n$ and $H'_n$ are simply the positive and negative parts of the
random walk $S_n$, and thus both are easily accessible -- refer to (\ref{esn}) below.

Since the study of large $N$ asymptotics is our primary objective,
we focus on the normal distribution for the demand model. Specifically, let
$$D_n=\mu_n+\sig Z_n,$$
where $Z_n$ are independent and identically distributed (i.i.d.) standard
normal random variables, and $\mu_n$ and $\sig$ are the mean and standard
deviation of $D_n$, respectively.
The common $\sig$ can be interpreted as a forecast error; hence, it is
independent of the periods.
Alternatively, we can allow a period-dependent $\sig_n$, as long as the
sum $\sum_{i=1}^n\sig_i \sim \sig\sqrt{n}$, i.e., it grows at a rate of
$\sqrt n$.
In this case, the above expression will still serve as a good approximation
due to the central limit theorem.

On the production side, we write, for $ n=1,\dots, N$,
\begin{eqnarray}
\label{supply}
x_1+\cdots +x_n=\mu_1+\cdots +\mu_n +\kappa\sig n^\alpha,
\end{eqnarray}
where $\kappa$ and $\al \in [0,1]$ are two policy parameters
(or decision variables).
Note that the (cumulative) production quantity in \eqref{supply} consists of
two parts: the first part matches the mean of the  (cumulative) demand; and
the second part  can be interpreted as {\it safety stock} to offset
demand variability.
The parameter $\alpha$ is often referred to as the ``safety factor''.
This factor addresses the issue of how the safety stock should match up with
the demand variability, where the latter grows on the order $\sqrt n$.
Recall that $\kappa$ can take on positive or negative values:
when it is positive, $\kappa\sig n^\alpha$ is truly the safety stock;
when $\kappa$ is negative, this part is a calculated under-production,
suitable for settings where the production cost is much higher relative to
the lost-sales penalty.
Since (\ref{supply}) implies that the production quantity for each period $n$
is given by
\begin{eqnarray}
\label{supply1}
x_n=\mu_n +\kappa\sig [n^\alpha- (n-1)^\alpha], \qquad n=1,\dots, N ,
\end{eqnarray}
then the cumulative shortage process $L_n$ can be estimated in different ways
depending on the value of $\al$.

\subsection{Useful Facts about the Standard Normal Variable}
\label{sec:calculations}

We summarize below some useful expressions and estimates associated
with the standard normal variable.  Let $Z$
denote the standard normal variate with density and distribution functions
denoted by $\phi (x)$ and $\Phi (x)$, respectively.
Define the ``shortfall function'' (or ``loss function'') as
\begin{eqnarray}
\label{g}
G(x):=\ex[(Z-x)^+]=\int_x^{+\infty}(z-x)\phi(z)dz =\phi(x)-x\bar\Phi(x) ,
\end{eqnarray}
where $\bar\Phi (x):=1-\Phi (x)$ and
the last equality follows directly from $\phi'(z)=-z\phi(z)$.

Now, we list properties of $G(x)$ that will be useful later for our purposes.
\begin{itemize}
\item[1.]
$G(x)$ is decreasing and convex in $x$ (since $(Z-x)^+$ is decreasing and convex in $x$).
\item[2.]
The part of $G(x)$ that has significant curvature is limited to the neighborhood of
the origin $x=0$. For large $x >0$,
we have $G(x)\approx 0$ and $G(-x)\approx x$.
\item[3.]
A direct derivation yields:
\begin{eqnarray}
\label{g_int}
2\int_a^b G(x)dx = \Phi(b) -\Phi(a)+bG(b)-aG(a).
\end{eqnarray}
\item[4.]
\begin{eqnarray}
\label{upperlower2}
0\le G(x)
\le \frac{\phi(x)}{x^2}, \qquad \forall x>0 .
\end{eqnarray}
\item[5.]
\begin{eqnarray}
\label{upperlower2asy}
\frac{x^2 G(x)}{\phi(x)}  \to 1 \qquad {\rm as}\; x\to +\infty.
\end{eqnarray}
\end{itemize}

To verify  \eqref{upperlower2} and \eqref{upperlower2asy}, observe that
\begin{eqnarray}
\label{upperlower}
\Big(1-\frac{1}{x^2}\Big)\frac{1}{x}\phi (x)
\le \bar\Phi(x)\le \frac{1}{x}\phi (x), \qquad x>0 ,
\end{eqnarray}
where the upper bound follows from
\begin{eqnarray*}
\int_x^\infty e^{-u^2/2} du
\le  \int_x^\infty \Big(1+\frac{1}{u^2}\Big)e^{-u^2/2} du
= \frac{1}{x}e^{-x^2/2}
\end{eqnarray*}
and the lower bound follows from
\begin{eqnarray*}
\int_x^\infty e^{-u^2/2} du
\ge  \int_x^\infty \Big(1-\frac{3}{u^4}\Big)e^{-u^2/2} du
= \Big(1-\frac{1}{x^2}\Big)\frac{1}{x}e^{-x^2/2}.
\end{eqnarray*}
The (second) inequality in \eqref{upperlower2} follows immediately
from the lower bound in (\ref{upperlower}).
Combining the upper and lower bounds in (\ref{upperlower}),
we have
$\bar\Phi(x)\approx \frac{1}{x}\phi (x)$ for large (positive) $x$.
This, along with l'H\^opital's rule, helps to verify \eqref{upperlower2asy}.

\subsection{Estimating $\mathbf{\ex[L_N]}$ for $\mathbf{\alpha=1}$}
\label{sec:one}


When $\alpha=1$ we have, from (\ref{supply1}),
$x_n=\mu_n +\kappa\sig$; thus, $D_n-x_n=-\kappa\sig +\sig Z_n$.
In this case, the partial sums $\{S_n\}$ constitute a random walk,
and $\ex[L_n]$ follows from standard results concerning the maximum of a random walk.
Specifically, by Spitzer's Identity
(see \cite{spitzer}; also, Ross \cite{ross}, Proposition 7.1.5):
\begin{eqnarray}
\label{sp}
\ex[L_N] = \sum_{n=1}^N \frac{1}{n} \ex[S_n^+],
\end{eqnarray}
where $S_n \, \eqd \, -\kappa\sig n +\sig \sqrt{n} Z$.
Hence,
\begin{eqnarray}
\label{esn}
\ex[S_n^+]=\sig\sqrt{n}\ex[(Z-\kappa\sqrt{n})^+]
=\sig\sqrt{n} G (\kappa\sqrt{n}).
\end{eqnarray}
Approximating the summation by integration, we have
\begin{eqnarray}
\label{sumint}
\sum_{n=1}^N \frac{1}{n} \ex[S_n^+]\approx
\int_0^N \frac{1}{t} \ex[S_t^+] dt
=\s \int_0^N \frac{1}{\sqrt{t}}G (\kappa\sqrt{t}) dt
=\frac{2\s}{\kappa} \int_0^{\kappa\sqrt{N}} G (u) du .
\end{eqnarray}
Therefore, combining the above with (\ref{g_int}), we obtain
\begin{eqnarray}
\label{eln}
\ex[L_N] \approx
\frac{\s}{\kappa} \Big(\Phi(\kappa\sqrt{N})-\frac{1}{2}\Big)+
\s\sqrt{N}G(\kappa\sqrt{N}) .
\end{eqnarray}

When $\kappa >0$, the above is increasing in $N$.
As $N\to\infty$, the second term on the right-hand side vanishes, and the
first term converges to $\frac{\s}{2\kappa}$.
When $\kappa <0$, the second term on the right-hand side of (\ref{eln})
quickly becomes $\sig |\kappa| N$ as $N$ grows, and the first term approaches
a constant $\frac{\s}{2|\kappa|}$.
These results then can be summarized in the following proposition.

\begin{pro}
\label{pro:al1}
{\rm
For $\alpha =1$ and large $N$, we have
\begin{eqnarray}
\label{al1}
\ex[L_N] \approx \frac{\s}{2\kappa}, \quad{\rm when}\quad \kappa >0;
\qquad
\ex[L_N] \approx \frac{\s}{2|\kappa|} +\sigma |\kappa| N, \quad{\rm when}\quad \kappa <0 .
\end{eqnarray}
}
\end{pro}

\medskip

Hence, setting $\alpha=1$, the expected shortage is bounded by a constant,
or grows linearly in $N$, corresponding to a positive or negative $\kappa$,
respectively.

\section{Estimating $\mathbf{\ex[L_N]}$ for $\mathbf{\alpha\in [0,1)}$}
\label{sec:bounds}

When $\alpha <1$, it is generally a difficult problem to calculate the exact
value of $\ex[L_N]$, and thus we develop bounds to reveal its asymptotic order
as $N\rightarrow \infty$.
We shall initially focus on $\kappa >0$ in the next two subsections, and
then summarize the corresponding results for $\kappa \le 0$ in the following
subsection.
Lastly, we consider a Brownian approximation for the case when $\alpha = \frac12$.

\subsection{Case of $\mathbf{\alpha \in (\frac12, 1)}$ and $\mathbf{\kappa > 0}$}

From (\ref{lost}), taking into account $S_0=0$, we have
$$L_n=\max_{0\le j\le n}\{S_j\}= \max_{1\le j\le n}\{S_j^+\}, \qquad n=1,\dots, N.$$
Jensen's inequality implies
\begin{eqnarray}
\label{elnlb}
\ex[L_N] \ge \max_{1\le n\le N}\{\ex[S_n^+]\},
\end{eqnarray}
where, similar to (\ref{esn}), we obtain
\begin{eqnarray}
\label{esn1}
\ex[S_n^+]=\sig\sqrt{n}\ex[(Z-\kappa {n}^{\alpha -1/2})^+]
=\sig\sqrt{n} G (\kappa {n}^{\alpha -1/2}) .
\end{eqnarray}

From (\ref{upperlower2}), we know that
\begin{eqnarray*}
\sig\sqrt{n} G (\kappa {n}^{\alpha -1/2}) \le
\frac{\sig\sqrt{n}\phi(\kappa {n}^{\alpha -1/2})}{\kappa^2 {n}^{2\alpha -1}}
:=\frac{\sig\sqrt{n}\phi(u)}{u^2},
\qquad{\rm where}\qquad u:=\kappa {n}^{\alpha -1/2} .
\end{eqnarray*}
Hence, $\ex[S_n^+]$ will decrease to $0$
(since $\phi(u)$ decreases to $0$ exponentially fast as $u\to+\infty$).
The values of $n$ that achieve the maximum in the lower bound of (\ref{elnlb})
must satisfy the following
optimality equation (whose solution need not be unique):
$$\hf n^{-1/2} G (y)= n^{1/2} \bar\Phi(y)\kappa (\alpha -1/2)n^{\alpha -3/2}.$$
This simplifies to
\begin{eqnarray}
\label{opty0}
G (y)= (2\alpha -1)y\bar\Phi(y)
\qquad{\rm or}\qquad
\phi (y)= 2\alpha y\bar\Phi(y) ,
\end{eqnarray}
and from (\ref{upperlower}) we know that $\phi(y)\approx y\bar\Phi(y)$ when
$y$ is large.
Since $\alpha >\hf$, the solution to the above equation exists.
Indeed, from (\ref{upperlower}), we have
$$y\bar\Phi(y)\ge \Big(1-\frac{1}{y^2}\Big)\phi(y), \qquad y>0 .$$
If we replace $y\bar\Phi(y)$ in (\ref{opty0}) by this lower bound, then the resulting solution
\begin{eqnarray}
\label{y0}
\bar y:= \Big(\frac{1}{2\alpha -1}\Big)^{\hf}
\end{eqnarray}
is an upper bound to the solution of the equations in (\ref{opty0}),
since $\phi (y)$ is decreasing in $y$ for $y >0$.
Therefore, in this case we obtain
\begin{eqnarray}
\label{elnlb2}
\ex[L_N]\ge \max_{0\le n\le N}\{\ex[S_n^+]\}\approx
\sig\Big(\frac{y}{\kappa}\Big)^{1/(2\alpha -1)} G(y),
\end{eqnarray}
where $y$ is the solution to (\ref{opty0}) and we know $y\in (0, \bar y)$.

For an upper bound, we have
$L_N =\max_{0\le n\le N} \{S_n\} \le \sum_{n=1}^N S_n^+$, which, along with
\eqref{esn1},
leads to
\begin{eqnarray}
\label{cal0}
\ex[L_N]  \le \sigma \sum_{n=1}^N \sqrt{n} G(\kappa n^{\al - \frac12})
\approx \sigma \int_0^N \sqrt{t} G(\kappa t^{\al - \frac12}) dt.
\end{eqnarray}
Using a transformation of variable $u= t^{\al -\frac12}$, together with
$G(x) \le \phi(x) /x^2$ from (\ref{upperlower2}), we obtain
\begin{eqnarray}
\label{cal}
\int_0^N \sqrt{t} G(\kappa t^{\al - \frac12}) dt
&=&
\frac{2}{2\al-1} \int_0^{N^{\al -\frac12}} u^{\frac{4-2\al}{2\al -1}} G(\kappa u) du
\nonumber\\
&\le& \frac{2}{(2\al-1)\kappa^2} \int_0^{N^{\al -\frac12}} u^{\frac{6-6\al}{2\al -1}} \phi(\kappa u) du
\nonumber\\
&=& \frac{2}{2\al-1} \kappa^{-\frac{3}{2\al-1}}  \int_0^{\kappa N^{\al -\frac12}} v^{\frac{6-6\al}{2\al -1}} \phi(v) dv
\nonumber\\
&\le & \frac{2}{2\al-1} \kappa^{-\frac{3}{2\al-1}}  \int_0^\infty v^{\frac{6-6\al}{2\al -1}} \phi(v) dv
:=C_\al <\infty .
\end{eqnarray}
Since $\frac{6-6\al}{2\al -1} >0$, as $\al\in (\hf, 1)$, we can summarize our
results for this case as follows.



\begin{pro}
\label{pro:al>0.5}
{\rm
For $1/2< \al <1$ and $\kappa >0$, we have\\
\hspace*{0.15in} (a) the lower bound:
\begin{eqnarray}
\label{lb>0.5}
\ex[L_N]\ge
\sig\Big(\frac{y^*_\al}{\kappa}\Big)^{1/(2\alpha -1)} G(y^*_\al), \quad{\rm for}\quad N\ge \bar y_\al:=(2\al -1)^{-1/2},
\end{eqnarray}
\hspace*{0.38in} where $y^*_\al>0$ is the solution to the equation $\phi (y)=2\al y\bar\Phi (y)$
(and hence, independent of $N$); and\\
\hspace*{0.15in} (b) the upper bound: $\ex[L_N]  \le \sig C_\al$, with $C_\al$
as specified in
(\ref{cal}).\\
Hence, in this case $\ex[L_N]$ is both upper- and lower-bounded by constants.
}
\end{pro}

\subsection{Case of $\mathbf{\alpha \in [0,\frac12]}$ and $\mathbf{\kappa > 0}$}


For $\alpha \in [0,\frac12]$, as $n\to\infty$, we have $G (\kappa {n}^{\alpha -1/2})\to G(0)$ or $G(\kappa)$
depending on whether $\al <\hf$ or $\al=\hf$.
Hence, from (\ref{elnlb}) and (\ref{esn1}), we obtain,
for sufficiently large $N$,
\begin{eqnarray}
\label{elnlb1}
\ex[L_N]\ge\sig G(\kappa)\sqrt{N}, \quad{\rm when}\; \al=\hf; \qquad
\ex[L_N]\ge\sig G(0)\sqrt{N}, \quad{\rm when} \; \al\in [0,\hf).
\end{eqnarray}



For an upper bound,
observe that the concave function $f(n)=n^{\alpha}$ is bounded from {\it below} by a linear function:
\begin{eqnarray}
\label{slowerrates}
n^\alpha \ge \frac{N^\alpha}{N} n:= \Theta n, \qquad {\rm for}\quad 0\le n\le N .
\end{eqnarray}
Defining $W_n:= -\kappa \sigma \Theta n + \sigma \sum_{k=1}^n Z_n$,
we know that $L_N \le \max_{n\le N} W_n$.
Meanwhile, $\max_{n\le N} W_n$ can be calculated through Spitzer's identity.
More specifically, we have
\begin{eqnarray*}
\ex[W_n^+] =\sig\sqrt{n} G (\kappa\Theta\sqrt{n}), \qquad 0\le n\le N,
\end{eqnarray*}
and
\begin{eqnarray*}
\sum_{n=1}^N \frac{1}{n} \ex[W_n^+]\approx
\int_0^N \frac{1}{t} \ex[W_t^+] dt
=\frac{2\s}{\kappa\Theta} \int_{0}^{\kappa\Theta \sqrt{N}} G(x)dx .
\end{eqnarray*}
Therefore, making use of \eqref{g_int}, we obtain
\begin{eqnarray}
\label{elnub0}
\ex[L_N] \le \sum_{n=1}^N \frac{1}{n} \ex[W_n^+]
\approx
\frac{\s}{\kappa\Theta} \Big(\Phi(\kappa\Theta \sqrt{N})-\hf\Big)
+{\s}\sqrt{N}G(\kappa\Theta \sqrt{N}).
\end{eqnarray}
As $N\to \infty$,
and according to the definition of $\Theta$ in (\ref{slowerrates}),
we have
$\Theta \sqrt{N} \to 0$ if $\alpha \in [0, \hf)$ and
$\Theta \sqrt{N} \to 1$ if $\al=\hf$.
The second term on the right-hand side of (\ref{elnub0}) is of order $\sqrt{N}$,
whereas the first term can be written as
$\frac{\s\sqrt{N}}{\kappa\Theta\sqrt{N}} [\Phi(\kappa\Theta \sqrt{N})-\hf]$.
Since $\lim_{y\to 0} (\Phi(y)-\hf)/y =\phi (0)$, we know this term is on
the order of $\sqrt{N}$ as well.
Hence,
\begin{eqnarray}
\label{elnub1}
\ex[L_N] \le
{\s}\sqrt{N}[ \phi (\kappa\Theta \sqrt{N})+G(\kappa\Theta \sqrt{N})] \le {\s}\sqrt{N}[ \phi (0)+G(0)]=2{\s}\phi (0)\sqrt{N}  .
\end{eqnarray}
These results allow us to conclude as follows.

\begin{pro}
\label{pro:al<0.5}
{\rm
In the case of $\alpha \in [0,\frac12]$ and $\kappa > 0$, $\ex[L_N]$ grows on
the order of $\sqrt{N}$.
In particular, we have the lower- and upper-bounds for $\ex[L_N]$  in  (\ref{elnlb1}) and (\ref{elnub1}).
}
\end{pro}

\subsection{Case of $\mathbf{\kappa \leq 0}$}

First consider the lower bounds.
For  $\al\in (\hf, 1)$ and $\kappa <0$, we have
\begin{eqnarray*}
\ex[S_n^+] = \sigma \sqrt{n} G(\kappa n^{\al -1/2})\approx \sig |\kappa| n^\al,
\end{eqnarray*}
since  $G(-x)\approx x$ when $x>0$ is large.
Hence, for $N$ sufficiently large,
\begin{eqnarray}
\label{lbneg}
\ex[L_N] \ge \max_{0\le n \le N} \{ \ex[S_n^+]\} \approx \sigma |\kappa| N^\al,
  \quad{\rm for} \; \al\in (\hf, 1) .
\end{eqnarray}
When $\al\in [0,\hf ]$ and $\kappa <0$, the lower bound for $\kappa>0$ in
(\ref{elnlb1}) remains valid; namely,
\begin{eqnarray}
\label{elnlbk<0}
\ex[L_N]\ge\sig G(\kappa)\sqrt{N}, \quad{\rm when}\; \al=\hf; \qquad
\ex[L_N]\ge\sig G(0)\sqrt{N}, \quad{\rm when} \; \al\in [0,\hf).
\end{eqnarray}

Next, we consider the upper bounds.
Define $Y_n := S_n +\kappa\sig n^\al$, for $n=0,1,\dots,N$.
Note that $Y_n$ is a  random walk with zero drift and $Y_0=0$.
Then,
\begin{eqnarray}
\label{ubk<0}
L_N&=&\max_{0 \le n\le N}\{ S_n\}
= \max_{0 \le n\le N}  \{Y_n - \kappa\sig n^\al \}
\nonumber\\
&\le&  \max_{0 \le n\le N}  \{Y_n\}  +  \max_{0 \le n\le N} \{- \kappa\sig n^\al \}
\nonumber\\
&=& \max_{0 \le n\le N}  \{Y_n\}  + (- \kappa\sig N^\al) ,
\end{eqnarray}
where the inequality follows from the subadditivity of the maximum operator,
and the last equality is due to $\kappa <0$.
Applying Spitzer's identity to the random walk $Y_n$,
similar to the derivation of $\ex[L_N]$ in (\ref{eln}) for the case of $\al=1$, we obtain
$$ \ex\Big[ \max_{0 \le n\le N}  \{Y_n \}\Big] =\ex\Big[ \max_{1 \le n\le N}  \{Y^+_n \} \Big]
= \sum_{n=1}^N \frac{\sig\sqrt{n}}{n} G(0)\approx \sig G(0) \int_0^N \frac{dt}{\sqrt{t}} =2\sig G(0) \sqrt{N}.
$$
Hence, combining this with \eqref{ubk<0}, renders
\begin{align}
\label{ubk<0a}
\ex[L_N] \le 2\sig G(0) \sqrt{N} +  |\kappa| \sig N^\al, \qquad \al \ge 0.
\end{align}
Our results for this case then can be summarized in the following proposition.

\begin{pro}
\label{pro:negk}
{\rm
For $\kappa<0$, $\ex[L_N]$ grows on the order of $\sig |\kappa| N^\al$ when
$\al\in(\hf, 1 )$;
it grows on the order of $\sig G(0) \sqrt{N} $ when $\al\in[0, \hf)$;
and it grows on the order of $\sig G(k) \sqrt{N} $ when $\al= \hf$.
In particular, the lower bound of $\ex[L_N]$ follows (\ref{lbneg}) and
(\ref{elnlbk<0}), with its upper bound following (\ref{ubk<0a}).
}
\end{pro}

\medskip

Finally, when $\kappa =0$, we have from (\ref{ubk<0}) and (\ref{ubk<0a})
$$\ex[L_N]=\ex\Big[ \max_{0 \le n\le N}  \{Y_n \}\Big] =2\sig G(0) \sqrt{N},$$
which is also consistent with setting $\kappa=0$ in the lower bound
(\ref{elnlbk<0}) and in the upper bound (\ref{ubk<0a}).

\subsection{Brownian approximation for case of $\mathbf{\alpha = \frac12}$}
\label{sec:hitting_time_approach}

The case of $\al= \frac12$ will turn out to produce the right trade-off between cost and service,
as we will see in the next section,
thus yielding the optimal order of production.
It is therefore appropriate to explore the possibility of a more accurate
estimation of $L_n$ in this particular case.

To this end, we study the continuous counterpart of $L_n$,
i.e., a Brownian motion with square-root drift.
More specifically, $\ex[L_N]$ can be approximated as
\begin{equation}
\ex[L_N] \;\; \approx \;\; \sqrt{\sigma}\ex\left[\sup_{1\le s\le t} \left\{B(s)-\frac{\kappa}{\sqrt{\sigma}}\sqrt{s}\right\}\right] ,
\label{eq:BrownianApproximation}
\end{equation}
where $B(t)$ is a standard Brownian motion.
Defining
\begin{align*}
\tau_x & := \inf \Big\{t > 1:B(t)=\frac{\kappa}{\sqrt{\sigma}}\sqrt{t}; B(1) = - x \Big\} 
\end{align*}
for any $x\ge 0$, we then have
\begin{eqnarray}
\ex\left[\sup_{1\le s\le t} \left\{B(s)-\frac{\kappa}{\sqrt{\sigma}}\sqrt{s}\right\}\right]&=& \int_0^\infty
\pr\left[\sup_{1\le s\le t} \left\{B(s)-\frac{\kappa}{\sqrt{\sigma}}\sqrt{s}\right\} \ge x \right] dx \nonumber \\
&=& \int_0^\infty \pr\left[\tau_x < t \right] dx.
\label{eq:BA:HT}
\end{eqnarray}

Now, define a new process $Y(u) := B(e^{2u})/e^u$.
One can readily verify, as shown in~\cite{Breiman}, that $Y(u)$ is an
Ornstein-Uhlenbeck process with $Y(0)=x$.
Further define
\begin{equation*}
T_{a,b} \;\; := \;\; \inf \{u\ge 0: Y(u) \ge b, Y(0) = a\} .
\end{equation*}
We note that, for any $t \ge 0$,
\begin{equation}
\pr[ \tau_x > e^{2t}] \;\; = \;\; \pr[ T_{-x,\kappa/\sqrt{\sigma}} >t] ,
\label{eq:BA:FPT}
\end{equation}
a fact that has been known since~\cite{Breiman}.
Our task is thus reduced to the evaluation of $T_{a,b}$.
When $b=0$, an explicit expression for $T_{a,b}$ has been derived;
see, e.g.,~\cite{RicciardiSato,leblanc,GJYor}.
More generally when $b\neq 0$ (as is of interest here), however,
only the Laplace transform of the density of $T_{a,b}$ is known.
Using various methods, the moments of this distribution also have
been calculated;
refer to, e.g.,~\cite{Siegert51,Sato78,CRS81,NRS85,RicciardiSato}.
In theory, the distribution function of $T_{a,b}$ can be recovered
from all of its moments.
In practice, the distribution function of $T_{a,b}$ can be
effectively approximated by its finite moments up to a certain
order.
Moreover, as shown in \cite{BertsimasPopescu}, the error bound of
such an approximation can be computed through numerical methods,
such as semi-definite programming.

We therefore propose to use the first $K$ moments $M_{a,b}^1, \ldots, M_{a,b}^K$
of $T_{a,b}$ where $a=-x$ and $b=\kappa/\sqrt{\sigma}$, in combination
with \eqref{eq:BrownianApproximation}, \eqref{eq:BA:HT} and \eqref{eq:BA:FPT},
to obtain a more accurate estimation of $L_n$ in the case of $\al= \frac12$.
The methods and results in~\cite{Siegert51,NRS85,RicciardiSato} can be
used to calculate the first $K$ moments
$M_{a,b}^1, \ldots, M_{a,b}^K$.
When $K \leq 3$, we can additionally exploit the closed-form tight
bounds in~\cite{BertsimasPopescu} to further improve our approximation
of the distribution function of $T_{a,b}$.
(Finding the best possible bounds for $K \geq 4$ is NP-hard~\cite{BertsimasPopescu}.)
As a specific instance of our proposed approach, we consider in more
detail the case of utilizing the first three moments of $T_{a,b}$
(i.e., $K=3$).
Our starting point is the set of closed-form expressions for
$M_{a,b}^1, M_{a,b}^2, M_{a,b}^3$ provided in~\cite{NRS85} and
parameterized by $a=-x$ and $b=\kappa/\sqrt{\sigma}$.
Then, letting $X$ denote a generic random variable that follows the
probability distribution function of $T_{-x,\kappa/\sqrt{\sigma}}$,
we next obtain closed-form tight upper bounds on both $\pr[X \le z]$
and $\pr[X\ge z]$ in terms of the expressions for
$M_{-x,\kappa/\sqrt{\sigma}}^1, M_{-x,\kappa/\sqrt{\sigma}}^2, M_{-x,\kappa/\sqrt{\sigma}}^3$.
More specifically, we exploit the following propostion adapted from
Theorem 3.3 in~\cite{BertsimasPopescu}.
\begin{pro}\label{prop:BertsimasPopescu}
For a constant $\delta > 0$ and a nonegative real random variable
$X \: \eqd \: T_{-x,\kappa/\sqrt{\sigma}}$ with first three moments
$M_1 = M_{-x,\kappa/\sqrt{\sigma}}^1$, $M_2 = M_{-x,\kappa/\sqrt{\sigma}}^2$
and $M_3 = M_{-x,\kappa/\sqrt{\sigma}}^3$, the functionals
$\pr[X > (1+\delta) M_1]$ and $\pr[X < (1-\delta) M_1]$ can be approximated as
\begin{equation*}
\pr[X > (1+\delta) M_1] \approx f_1(C_M^2,D_M^2, \delta) \qquad \mbox{and} \qquad \pr[X < (1-\delta) M_1] \approx f_2(C_M^2,D_M^2, \delta) ,
\end{equation*}
respectively, where
\begin{align*}
f_1(C_M^2,D_M^2, \delta)& = \left\{
  \begin{array}{cc} \min \left(\frac{C_M^2}{C_M^2+\delta^2}, \frac{1}{1+\delta} \cdot \frac{D_M^2}{D_M^2 + (C_M^2-\delta)^2}\right), & \qquad \delta > C_M^2 , \\
  \frac{1}{1+\delta} \cdot \frac{D_M^2+(1+\delta)(C_M^2-\delta)}{D_M^2 + (1+C_M^2)(C_M^2-\delta)}, & \qquad \delta \leq C_M^2 , \end{array} \right. \\
f_2(C_M^2,D_M^2, \delta)& = 1- \frac{(C_M^2+\delta)^3}{(D_M^2 + (C_M^2+1)(C_M^2+\delta))(D_M^2 + (C_M^2+\delta)^2)} , \qquad \delta < 1 , \\
C_M^2 & = \frac{M_2-M_1^2}{M_1^2} , \\
D_M^2 & = \frac{M_1M_3-M_2^2}{M_1^4} .
\end{align*}
\end{pro}
Lastly, we combine Proposition~\ref{prop:BertsimasPopescu} and the expressions for
$M_{-x,\kappa/\sqrt{\sigma}}^1, M_{-x,\kappa/\sqrt{\sigma}}^2, M_{-x,\kappa/\sqrt{\sigma}}^3$
in \cite{NRS85} together with \eqref{eq:BrownianApproximation}, \eqref{eq:BA:HT},
\eqref{eq:BA:FPT} to obtain the more accurate approximation of $L_n$ in the case
of $\al= \frac12$, as desired.

\section{Cost Minimization and Asymptotically Optimal Solutions}
\label{sec:optimality}

In this section we consider instances of two cost minimization models for the lost-sale system of interest.
For the first model, the total cost consists of both the production cost and lost-sale penalty.
To be more precise,
we minimize $c\sum_{n=1}^N x_n + p \ex[L_N]$ where $c$ and $p$ are the per-unit production cost and lost-sales penalty, respectively.
From \eqref{supply},
ignoring the constant term $c(\mu_1+\cdots +\mu_N)$, we have
\begin{eqnarray}
\label{prodlost}
\min_{\alpha,\kappa} \quad c\sig\kappa N^\al + p \ex[L_N] .
\end{eqnarray}
For the second model, we incorporate the inventory (holding) cost in each period.
Specifically, we have
\begin{eqnarray}
\label{prodlostinv}
\min_{\alpha,\kappa} \quad c\sig\kappa N^\al + p \ex[L_N] + h \ex[H_N] ,
\end{eqnarray}
where $h$ is the per-unit inventory holding cost.
Recall from Proposition \ref{pro:lost} that the inventory (at the end) of
period $n$ is given by $H_n=L_n-S_n$, and therefore
$\ex[H_N]= \ex[L_N]+\kappa \sigma N^\alpha$.

Our first step is to find the optimal value of $\al$ in these optimization problems.
For this purpose,
let us start with a summary of the results we have derived so far concerning the asymptotics of
the expected loss sales in the limit as $N\to\infty$.
\begin{itemize}
\item[(i)]
When $\alpha=1$, $\ex[L_N]$ is bounded by a constant
$\frac{\sig}{2\kappa}$ if $\kappa >0$; whereas it grows linearly in $N$
if $\kappa <0$.
\item[(ii)]
When $\alpha\in (\hf, 1)$ and $\kappa >0$, $\ex[L_N]$ is bounded from above
and below by constants that are independent of $N$ (but dependent on $\al$
and $\kappa$); whereas it grows on the order of $N^\al$ if $\kappa <0$.
\item[(iii)]
When $\alpha\in [0,\hf ]$, $\ex[L_N]$ grows on the order of $\sqrt {N}$,
regardless of whether $\kappa >0$ or $\kappa \leq 0$.
\end{itemize}

From \eqref{prodlost},
if $\al \in (\hf,1]$ and $\kappa >0$, then the second term in the objective will be bounded by a constant independent of $N$;
hence, for sufficiently large $N$, the objective is of order $N^\al$, i.e., the same order as the first term.
On the other hand,
if $\al \in (\hf,1]$ and $\kappa <0$, then the first term will decrease while the second term will increase, both on the order of $|\kappa| N^\al$;
hence, if $c> p$, then $\kappa =-\infty$ minimizes the objective value, whereas if $c\le  p$, then $\kappa \to 0$ minimizes the objective value.

Therefore, the only non-trivial solution of the optimization problem is to set
$\al\le \hf$, corresponding to the case summarized in (iii) above.
Here, the first term of the objective function is of order $N^\alpha$, which
is dominated by the second term of order $\sqrt {N}$ (except when $\alpha=\hf$).
Consequently, the objective value will also be of order $\sqrt{N}$, as is
further confirmed by replacing $\ex[L_N]$ with its lower bound via Jensen's
inequality (which applies in all cases):
$\ex[L_N] \ge \ex[\max_{1\le n\le N} S_n^+]\ge \ex[S_N^+]$.
This leads to the following minimization problem, which is a lower bound of
the original problem:
\begin{eqnarray}
\label{eqn:production_LB_1}
\min_z \quad cz + p \ex\left[\left( \sig\sum_{n=1}^NZ_n -z\right)^+\right] = cz+p\ex \big[(\sig \sqrt{N} Z- z\big)^+].
\end{eqnarray}
To allow for any $\al\in[0,1]$ and any $\kappa$, we write
$$z=\sig\kappa N^\al =\sig\sqrt{N}\kappa N^{\al-\hf}:= \sig\sqrt{N} y,$$
and then the optimal solution to minimizing $cy+p \ex[(Z-y)^+]$ is given by
\begin{eqnarray}
\label{ylb}
y^*=\bar\Phi^{-1}\Big(\frac{c}{p}\Big)=\Phi^{-1}\Big(\frac{p-c}{p}\Big).
\end{eqnarray}
With $y^*$ being a constant, independent of $N$, we must have that $\al=\hf$,
and hence $\kappa =y^*$.
The corresponding objective value then can be expressed as
\begin{eqnarray}
\label{yoptval}
\sigma \sqrt{N}  [c y^* +pG (y^*)]=\sigma \sqrt{N}  [c y^* +p\phi (y^*)-py^*\bar\Phi(y^*)]= \sigma \sqrt{N} p\phi(y^*) .
\end{eqnarray}

This confirms our previous statements; namely, the minimal overall cost cannot be
lower than order $\sqrt{N}$.
Moreover, we shall henceforth assume that $p\ge c$, unless noted otherwise,
because if $p<c$ then the lower-bound solution in (\ref{ylb}) already indicates
what will happen: make $y$ as small (negative) as possible; whereas the
objective value is lower bounded by the expression in (\ref{yoptval}).
We therefore only need to consider the case of $\al =\frac12$.
Furthermore, upon careful examination of the objective of the second optimization
model, we can reach the same conclusion.

In the remainder of this section, we start by considering asymptotically optimal solutions of
the above optimization models using two different approaches to obtain the desired solutions.
We then turn to incorporate inventory costs in our optimization models.
Finally, we consider the asymptotic equivalence of lost-sales and backorder models
under certain conditions.

\subsection{Asymptotically Optimal Solutions}

With our focus on $\alpha = \frac12$, we first obtain an asymptotically optimal solution through a Brownian approximation
of the objective function for large $N$.
We then obtain optimal solutions through an analysis of the asymptotic behavior of our upper and lower bounds of $\ex[L_N]$.

\subsubsection{Brownian Approximation}
\label{sec:Brownian_Approx}

Recall that the optimization problem for the production loss trade-off, when $\alpha$ is fixed to be $\frac12$, has the form:
\begin{eqnarray}
\label{prodlost_recall}
\min_{\kappa} \quad c\sig\kappa N^{1/2} + p \ex[L_N] .
\end{eqnarray}
Taking a closer look at the objective function, we conclude
\begin{align*}
c\sig\kappa N^{1/2} + p \ex[L_N]   & = \left(c\sig\kappa + p \frac{\ex[L_N]}{ N^{1/2}}\right)N^{1/2} .
\end{align*}
Moreover,
\begin{align*}
 \frac{\ex[L_N]}{ N^{1/2}} =\ex\left[\max_{n\le N}\frac{S_n}{N^{1/2}}\right]= \ex\left[\max_{n\le N}\frac{\sum_{i=1}^n\sig Z_i-\kappa \sig \sqrt{n} }{N^{1/2}}\right]
\end{align*}
where $Z_i$ are i.i.d.\ copies of a random variable following the standard normal distribution. 
This leads to our next result.
 
\begin{thm}
Let $\rho(\kappa)= \sup_{t\le 1 } B_t -\kappa \sqrt{t}$, where $B_t $ is a standard Brown motion.  Then,
\begin{align*}
\lim_{N\rightarrow \infty}  \ex\left[\max_{n\le N}\frac{\sum_{i=1}^nZ_i-\kappa  \sqrt{n} }{N^{1/2}}\right] \rightarrow \rho(\kappa) .
\end{align*}
\end{thm} 
\proof{Proof.}
We first obtain
\begin{align*}
\max_{n\le N}\frac{\sum_{i=1}^nZ_i-\kappa \sqrt{n} }{N^{1/2}}& = \max_{t \in \{ \frac{n}{N}, n=1,2,\ldots, N\}} \frac{\sum_{i=1}^{Nt}Z_i-\kappa \sqrt{Nt} }{N^{1/2}}\\ & = \sup_{t \in (0,1]} \frac{\sum_{i=1}^{\lfloor Nt\rfloor}Z_i-\kappa \sqrt{\lfloor Nt \rfloor} }{N^{1/2}} .
\end{align*}
From the functional central limit theorem, for the summation of random variables (see e.g.,~\cite{chenyaobook}), we can conclude that
$$\frac{\sum_{i=1}^{\lfloor Nt\rfloor}Z_i-\kappa  \sqrt{\lfloor Nt \rfloor} }{N^{1/2}} \Rightarrow B_t -\sqrt{t}.$$
The desired convergence then follows from the continuous mapping theorem.
\Halmos\endproof

The above theorem implies that, for any $\epsilon>0$ and when $N$ is large enough, we have
\begin{align*}
\Big| \frac{\ex[L_N]}{ N^{1/2}}-\rho(\kappa)\Big| < \epsilon
\end{align*}
uniformly on compact sets in $\kappa$.
Although there is no explicit formula for $\rho(\kappa)$, computational methods that make use of the related results in \cite{APP,Breiman,RicciardiSato,Sato78}
together with Monte-Carlo simulation can be employed;
in particular, these methods can be used to obtain accurate estimation of the value of $\rho(\kappa)$ and its derivative in order to solve the Brownian version of our optimization problem:
\begin{equation}
\min_{\kappa} \quad  c\sig\kappa N^{1/2} + p \sig \rho(\kappa).
\label{BrownOpt}
\end{equation}
Let $\kappa^*$ denote an optimum of the problem \eqref{BrownOpt}; of course, $\kappa^*$ corresponds to a specific production plan.
We next show that this production plan is asymptotically optimal, defined as follows.
\begin{defn}
A production plan is {\it asymptotically optimal} if
\begin{align*}
 \lim \sup_{N\rightarrow \infty} \frac{J_N(\kappa)}{J_N^*} \le 1. 
\end{align*} 
\end{defn}

\begin{thm}
The optimal solution of \eqref{BrownOpt}, $\kappa^*$, is asymptotically optimal for the problem \eqref{prodlost_recall}.
\end{thm}
\proof{Proof.}
Let $\kappa^*(N)$ denote an optimum of problem \eqref{prodlost_recall} for a given $N$.
Then, for any $N$, we have
\begin{align*}
 \frac{J_N(\kappa)}{J_N^*} &= \frac{c\sig\kappa^* N^{1/2} + p \ex[L_N(\kappa^*)] } {c\sig\kappa^*(N) N^{1/2} + p \ex[L_N(\kappa^*(N))] }\\  & =
\frac{c\sig\kappa^* N^{1/2} + p \ex[L_N(\kappa^*)] - \left(c\sig\kappa^*(N) N^{1/2} + p \ex[L_N(\kappa^*(N))]\right)+ \left(c\sig\kappa^*(N) N^{1/2} + p \ex[L_N(\kappa^*(N))]\right)} {c\sig\kappa^*(N) N^{1/2} + p \ex[L_N(\kappa^*(N))] }
\\ & \le \frac{c\sig\kappa^* N^{1/2} + p \ex[L_N(\kappa^*)] - \left(c\sig\kappa^*(N) N^{1/2} + p \ex[L_N(\kappa^*(N))]\right)+ } {c\sig\kappa^*(N) N^{1/2} + p \ex[L_N(\kappa^*(N))] } + 1.
\end{align*}
We further know that the first term can be bounded by a constant, and thus goes to $0$ in the limit at $N \rightarrow \infty$.
The desired result of asymptotic optimality follows.
\Halmos\endproof

It is evident that the asymptotically optimal result is not restricted to the production optimization problem, as a similar conclusion can be also reached for the inventory optimization problem.

\subsubsection{Lower- and Upper-Bound Surrogate Problems}
\label{sec:prodlost}

%

We now turn to consider an alternative approach for obtaining optimal solutions through our upper and lower bounds of $\ex[L_N]$.
For finite $N>0$, we have already analyzed the lower bound surrogate problem in (\ref{yoptval}).
We next seek to minimize the upper bound with respect to the choice of $\kappa$.
From the upper bound for the case $\al=\hf$ and $\kappa >0$ in terms of
\eqref{slowerrates} and \eqref{elnub0},
upon replacing $c\sigma\sqrt{N}$ and $p\sigma\sqrt{N}$ with $c$ and $p$,
respectively, we want to solve the following optimization problem:
\begin{eqnarray}
\label{ubobj}
\min_{\kappa \ge 0} \; c\kappa + p \left( \frac{1}{\kappa} \Big[\hf - \bar\Phi(\kappa)\Big]+G(\kappa)\right).
\end{eqnarray}
The optimality condition is then given by
$$c+\frac{p}{\kappa^2} \Big(\kappa \phi(\kappa) +\bar\Phi(\kappa) -\hf \Big) = p \bar\Phi(\kappa), $$
which simplifies to
\begin{eqnarray}
\label{uboptsol}
c\kappa + p G(\kappa) = \frac{p}{\kappa} \Big(\hf - \bar\Phi(\kappa)\Big) .
\end{eqnarray}
Denoting the solution of \eqref{uboptsol} by $\kappa_u$ and observing
$\kappa_u\ge 0$ implies that $\bar\Phi(\kappa_u)\le \hf$.
Further denote the objective values of the lower- and upper-bound problems
in (\ref{yoptval}) and (\ref{ubobj}) by $V^\ell$ and $V^u$, respectively.
We write $\kappa_\ell:=y^*$ following (\ref{ylb}), which yields
\begin{eqnarray}
\label{kl}
\kappa_\ell :=\bar\Phi^{-1}\Big(\frac{c}{p}\Big)=\Phi^{-1}\Big(\frac{p-c}{p}\Big).
\end{eqnarray}
Then, from (\ref{yoptval}), (\ref{ubobj}) and \eqref{uboptsol}, we have
\begin{eqnarray}
\label{vuvell}
\frac{V^u}{V^\ell} = \frac{2[ c\kappa_u +pG(\kappa_u)] }{c\kappa_\ell +pG(\kappa_\ell)}  \ge 2,
\end{eqnarray}
since $\kappa_\ell$ is the minimizer of $cy+pG(y)$.

To go in the opposite direction,
from (\ref{yoptval}), (\ref{ubobj}) and (\ref{uboptsol}), we can also write
the above ratio as
$$\frac{V^u}{V^\ell} = \frac{2 [\hf - \bar\Phi(\kappa_u)]}{\kappa_u \phi(\kappa_\ell)}.$$
Observe that $\frac{1}{x} [\hf - \bar\Phi(x)]$ is decreasing in $x\ge 0$
(and at $x=0$, the limit is $\phi (0)$, via l'H\^opital's rule), since its
derivative is $\frac{1}{x^2} [x\phi(x) + \bar\Phi(x)-\hf]$.
It is then straightforward to verify that $x\phi(x) + \bar\Phi(x)-\hf\le 0$
for $x\ge 0$.
We therefore have
$$\phi (0)\ge \frac{1}{\kappa_u} \Big(\hf - \bar\Phi(\kappa_u)\Big) ,$$
and hence
\begin{eqnarray}
\label{vuvell1}
\frac{V^u}{V^\ell} \le \frac{2\phi(0)}{\phi(\kappa_\ell)} .
\end{eqnarray}
In fact, we also know that $\kappa_\ell\le\kappa_u$, which follows from
$$\frac{p}{\kappa_\ell} \Big(\hf - \bar\Phi(\kappa_\ell)\Big)\ge \frac{p}{\kappa_u} \Big(\hf - \bar\Phi(\kappa_u)\Big)$$
(since $\kappa_u$ is the minimizer of the right-hand side), together with
the fact  justified above that $\frac{1}{x} [\hf - \bar\Phi(x)]$ is decreasing in $x$.
When $p=2c$, we have $\kappa_\ell =0$; refer to (\ref{kl}).
Then, we will have $\kappa_u=0$ as well; refer to (\ref{uboptsol}).
These results render, for this case, $V^u=2V^\ell$.

In summary, we have established the following proposition.
\begin{pro}
\label{pro:lostbounds}
{\rm
Let  $V^*$ denote the optimal value of the original problem in (\ref{prodlost}).
Let $V^\ell$ and $V^u$ denote the objective values of the lower- and
upper-bound problems in (\ref{yoptval}) and (\ref{ubobj}), respectively.
Then,
$$  \frac{V^u}{V^\ell}= \frac{2[ c\kappa_u +pG(\kappa_u)] }{c\kappa_\ell +pG(\kappa_\ell)}
=\frac{2 [\hf - \bar\Phi(\kappa_u)]}{\kappa_u \phi(\kappa_\ell)} $$
where $\kappa_\ell=y^*$ is the lower-bound solution following (\ref{ylb}),
$\kappa_u$ is the upper-bound solution to the optimality equation in
(\ref{uboptsol}), and $\kappa_\ell\le\kappa_u$.
Consequently, we have
$$2\le \frac{V^u}{V^\ell} \le \frac{2\phi(0)}{\phi(\kappa_\ell)}$$
and
$$\frac{V^\ell}{V^u}\le \frac{V^\ell}{V^*}\le \frac{V^u}{V^*} \le \frac{V^u}{V^\ell}.$$
}
\end{pro}

\medskip

Based on the foregoing analysis, we propose to use $\al=\hf$ and $\kappa=\kappa_u$,
the solution to the upper-bound problem $V^u$ (for $\kappa >0$), for our
approximate solution to the original optimization problem in (\ref{prodlost}).
When $\kappa <0$,  we can still use the solution to $V^u$ as a heuristic, simply
because the upper- or lower-bounds in this case will not lead to meaningful solutions.
(For instance, the upper-bound problem is $\min_{\kappa\le 0} (c-p)\kappa$,
and thus the solution is: $\kappa^*=0$ if $c\le p$, and $\kappa^* =-\infty$
if $c>p$.)
All we need is to remove the constraint $\kappa \ge 0$ from the minimization
problem in (\ref{ubobj}), and then the optimality equation in (\ref{uboptsol})
still applies (with the right-hand side {\it increasing} in $\kappa$ for
$\kappa\le 0$).
However, it is no longer the case that $V^u$ provides an upper bound.
In fact, when $\kappa<0$, $V^u$ is a {\it lower} bound of $V^*$: the
constructed linear lower bound to the concave function $\sqrt n$, along with a
negative $\kappa$, implies subtracting {\it less} from the production costs,
which in turn must result in lower lost-sales penalties.
Thus, we have $V^u/V^* \le 1$, and the performance of $V^u$, as an
approximation to $V^*$, should be comparable to, if not better than, its
performance in the case of $\kappa >0$.


\subsection{Incorporating Inventory Cost}
\label{sec:inventory_cost}

Now, we turn to incorporate inventory costs in our analysis.
Referring back to \eqref{prodlostinv},
let $h/N$ be the holding cost per unit of inventory that is held for one period,
with $h>0$ a constant parameter.
The reason for including $N$ in the denominator can be explained as follows.
The application context we have in mind is a fast-moving production-inventory
system, where the planning horizon is of a fixed length, say one unit, and $h$
is the cost to hold one unit of inventory over this entire horizon.
By dividing the horizon into $N$ segments (periods), we are effectively using
$N$ as a scaling parameter on the demand.
Hence, as previously noted, instead of demand $D$ over the entire planning horizon,
we have been considering $D_n$ for each period $n=1,\dots, N$, with $D_n\eqd D$ and
independent of $n$.
Another way to motivate this approach is due to the fact that, if $h$ is not
scaled by $N$, then the inventory holding cost, because of its cumulative
nature over time, will dominate all other costs when $N$ grows.

Therefore, using $\al=\hf$, the inventory cost is given by
$$\frac{h}{N}\sum_{n=1}^N \sig\sqrt{n} [G(\kappa) +\kappa]\approx
\frac{2}{3} h  \sig\sqrt{N} [G(\kappa) +\kappa]$$
for the lower bound;
and similarly, for the upper bound, we have
$$\frac{h}{N}\sum_{n=1}^N \sig\sqrt{n} \left(\frac{1}{\kappa} \Big[\hf-\bar\Phi(\kappa)\Big]
+G(\kappa) +\kappa\right)\approx
\frac{2}{3} h  \sig\sqrt{N}  \left(\frac{1}{\kappa} \Big[\hf-\bar\Phi(\kappa)\Big]
+G(\kappa) +\kappa\right);$$
refer to \eqref{elnlb1} and \eqref{elnub0}.
Upon replacing $\frac{2}{3} h$ by $h$, the inventory costs derived above are
easily incorporated into both the Brownian approximation and the lower- and upper-bound
surrogate problems of the previous subsections:
simply replace $c$ by $c+h$, and $p$ by $p+h$.

\subsection{Asymptotic Equivalence of Lost-Sales and Backorder Models}
\label{sec:backorder}

Next, we consider another type of asymptotic results.
Suppose the penalty $p$ is large, say $p\to\infty$, whereas $N$ is fixed.
We first ask the question of whether it is possible to keep the expected
penalty cost $p\ex[L_N]$ finite (for a given $N$).

To this end,
let us consider the lower-bound solution with $\al=\hf$ and
$\kappa=\kappa_\ell=y^*$;
for simplicity, we shall write $y$ instead of $y^*$.
Recall from \eqref{ylb} that $y =\bar\Phi^{-1}(\frac{c+h}{p+h})$ when inventory
cost is included into the model, and thus $p\to\infty$ is equivalent to
$y\to\infty$.
Then, for any $n$, we have
$$p\ex[S_n^+] =\sig\sqrt{n} pG(y)=\sig\sqrt{n} \frac{cG(y)}{\bar\Phi(y)},$$
for which applying l'H\^opital's rule yields
$$\lim_{y\to\infty}\frac{G(y)}{\bar\Phi (y)}=\lim_{y\to\infty}\frac{\bar\Phi (y)}{\phi(y)}=\lim_{y\to\infty}\frac{1}{y}\to 0.$$
Namely, $p\ex[S_n^+]\to 0$ for every $n$.

In fact, with $S^+_n(y):=\sig\sqrt{n}(Z-y)^+$, this strategy also leads to
$\ex[S_n^+(y)]=\sig\sqrt{n} G(y)\to 0$ as $y\to\infty$ for every $n$.
Because $S^+_n(y)$ is decreasing in $y$, monotone convergence implies that
$S^+_n(y) \to 0$ as $y\to\infty$, for every $n=1,\dots, N$.
Hence, we must have $L_N(y)=\max_{1\le n\le N} S_n^+ (y)\to 0$ as $y\to\infty$,
and thus $\ex[L_N(y)]\to 0$ as $y\to\infty$.
Since the optimal strategy (to the original problem with inventory cost
included) can do no worse, we must have
\begin{eqnarray*}
&& \min_\kappa \Big(\sig\sqrt{N} (c+h)\kappa + (p+h)\ex[L_N (\kappa) ]\Big)
\nonumber\\
&\le& \sig\sqrt{N} (c+h)y + (p+h)\ex[L_N (y) ]
= \sig\sqrt{N} (c+h)y+o(y),
\end{eqnarray*}
where $y=\bar\Phi^{-1}(\frac{c+h}{p+h})$ as specified above.

On the other hand, based on the results established earlier, we know that the
right-hand side above is also a lower bound of the objective function, because
(with $\kappa^*$ denoting the minimizer of the first expression)
\begin{eqnarray*}
&& \min_\kappa  \Big(\sig\sqrt{N}  (c+h)\kappa + (p+h)\ex[L_N (\kappa) ]\Big) \\
& \ge& \sig\sqrt{N} (c+h)\kappa^* + (p+h)\ex[S^+_N (\kappa^*) ]\\
& \ge& \sig\sqrt{N} (c+h)y + (p+h)\ex[S^+_N (y) ]
= \sig\sqrt{N}  (c+h)y+o(y),
\end{eqnarray*}
where the second inequality follows from $y$ being the minimizer of
the lower-bound problem.
We therefore have
\begin{eqnarray}
\label{lsub}
\min_\kappa \Big(\sig\sqrt{N} (c+h)\kappa + (p+h)\ex[L_N (\kappa) ]\Big)
= \sig\sqrt{N} (c+h)y+o(y).
\end{eqnarray}

Next, let us consider the backorder model.
Recall that $B_n=S_n^+$ is the number of backlogged units of demand and
that $H'_n=S_n^-$ is the inventory, both at the end of period $n$;
refer to Proposition \ref{pro:lost}.
Similar to the lost-sales model in the previous subsection, let $\frac{b}{N}$
and $\frac{h'}{N}$ be the per-unit backlog penalty and inventory holding costs
in each period, respectively.
Then, the (original) optimization problem can be expressed as
\begin{eqnarray}
\label{bk}
\min_\kappa \left( c\kappa + \frac{b}{N}\sum_{n=1}^N \ex[S^+_n (\kappa) ]
+\frac{h'}{N}\sum_{n=1}^N \ex[S_n^- (\kappa) ] \right)
= \sig\sqrt{N} \,  \min_\kappa \left[ (c+h')\kappa + (b+h') G (\kappa) \right] ,
\end{eqnarray}
where $\frac{2}{3} b$ and $\frac{2}{3} h'$ are replaced by $b$ and $h'$
(analogous to the approach taken with $h$ in the previous subsection),
respectively, and taking into account
$S_n^-=S_n^+-S_n=S_n^+ +\sig\sqrt{n}\kappa$.
Let $y'$ denote the optimal solution to \eqref{bk}, clearly rendering
$y'=\bar\Phi^{-1}(\frac{c+h'}{b+h'})$.
Suppose we set $b=p$ and $h'=h$, and let $p\to\infty$.
Then, $y'=y$ and the left-hand side of (\ref{bk}) has the following
asymptotics for large $y$:
\begin{eqnarray}
\label{bklimit}
\sig\sqrt{N}  \min_\kappa \Big( c\kappa + b \ex[S^+_n (\kappa) ] +h' [S_n^- (\kappa) ] \Big)
= (c+h)y' + o(y') =(c+h)y + o(y).
\end{eqnarray}

Upon comparing (\ref{lsub}) and (\ref{bklimit}), we have established the following result.

\begin{pro}
\label{pro:equiv}
{\rm
As $p \rightarrow \infty$, the objective values of the lost-sales and
backorder models become equivalent in that their ratio approaches $1$,
provided the lost-sales and backlog penalty costs are both equal to $p$,
and the purchasing and holding costs are also respectively equal in the
two models.
}
\end{pro}

\section{Conclusions}
\label{sec:conclusions}

A general class of high-volume, fast-moving production-inventory problems are studied.
Although we focus on the case in which the production is set to be the average demand
supplemented by a certain amount of safety stock, our assumptions on this safety stock
is quite general.
In the asymptotic regime, we identify the right order of the safety stock by studying
random walks with power drifts.
We also derive bounds and approximations for functionals of key performance measures
so that the production planning can be optimized under different settings.
Finally, our analysis provides another means and settings for establishing the asymptotic
equivalence of lost sale and backorder models, as observed by several authors using
different means and settings.

\end{document}